\documentclass[11pt]{article}
\usepackage{amsmath, amsthm}

\usepackage{amsmath}
\usepackage{color}
\usepackage{amssymb}
\usepackage[stable]{footmisc}

\theoremstyle{plain}
\newtheorem{theorem}{\bf Theorem}[section]

\theoremstyle{definition}

\newtheorem{definition}[theorem]{\bf Definition}
\newtheorem{remark}[theorem]{\bf Remark}

\newcommand{\N}{\mathbb N}
\newcommand{\Z}{\mathbb Z}
\newcommand{\R}{\mathbb R}

\numberwithin{equation}{section}

\begin{document}

\title{Minimal $P$-symmetric periodic solutions of nonlinear Hamiltonian systems}

\author{
Shanshan Tang\thanks{Partially supported by initial Scientific Research Fund of Zhejiang Gongshang University. E-mail: ss.tang@foxmail.com}\\[1ex]
 School of Statistics and Mathematics, Zhejiang Gongshang University\\
 Hangzhou 310018, P.R. China
}

\date{}

\maketitle

\begin{abstract}
{In this paper some existence results for the minimal P-symmetric periodic solutions are proved for first order autonomous Hamiltonian systems when the Hamiltonian function is superquadratic, asymptotically linear and subquadratic. These are done by using critical points theory, Galerkin approximation procedure, Maslov $P$-index theory and its iteration inequalities.}
\end{abstract}



\section{Introduction and main results}
\bigskip

\ \ \ \ We study the $P$-boundary problem of first order autonomous Hamiltonian systems:
\begin{equation}\label{1}
\left\{
\begin{array}{l}
\dot{x}=JH^{\prime}(x),\ \forall x\in\R^{2n},\\
\!x(\tau) = Px(0),
\end{array}\right.
\end{equation}
where $\tau>0$, $P\in Sp(2n)$, $H\in C^{2}(\R^{2n}, \R)$ and $H(Px) = H(x)$, $\forall x \in \R^{2n}$. $H^{\prime}(x)$ denote its gradient, $J=\left( \begin{array}{cc}
0 \ \ & -I_{n}\\
I_{n} & 0
\end{array} \right)$
is the standard symplectic matrix, $I_{n}$ is the identity matrix on $\R^{n}$ and $n$ is the positive integer.

A solution $(\tau, x)$ of the problem (\ref{1}) is called {\it $P$-solution} of the Hamiltonian systems. It is a kind of generalized periodic solution of Hamiltonian systems. The problem (\ref{1}) has relation with the the closed geodesics on Riemannian manifold (cf.\cite{HuSun2}) and symmetric periodic solution or the quasi-periodic solution problem (cf.\cite{HuSun}). In addition, C. Liu in \cite{LC5} transformed some periodic boundary problem for asymptotically linear delay differential systems and some asymptotically linear delay Hamiltonian systems to $P$-boundary problems of Hamiltonian systems as above, we also refer \cite{CAMR,FDTS, HuSun1,HuWang} and references therein for the background of $P$-boundary problems in $N$-body problems.

Suppose $P$ satisfies $P^{k}=I$, here $k$ is assumed to be the smallest positive integer such that $P^{k}=I$ (this condition for $P$ is called $(P)_{k}$ condition in the sequel), so the {\it $P$-solution} $(\tau, x)$ can be extended as a $k\tau$-periodic solution $(k\tau, x^{k})$. We say that a $T$-periodic solution $(T, x)$ of a Hamiltonian system in (\ref{1}) is $P$-symmetric if $x(\frac{T}{k})=Px(0)$. $T$ is the {\it $P$-symmetric} period of $x$. We define $T$ be the minimal {\it $P$-symmetric} period of $x$ if $T=\min\{\lambda>0\mid x(t+\frac{\lambda}{k})=Px(t), \forall t\in \R\}$. Note that $T$ might not be the minimal period of $x$ although it is the minimal {\it $P$-symmetric} period of $x$.

In recent years, Maslov P-index theory was developed to study the existence and multiplicity of {\it $P$-solutions} (cf.\cite{dong,donglong,LC,LC5}), specially, the corresponding iteration theory was built to estimate the minimality of the period of {\it $P$-solution} (i.e., the minimal {\it $P$-symmetric} period) (cf.\cite{LC2,liutss2}) and look for geometrically distinct $P$-solutions (i.e., subharmonic $P$-solutions) (cf.\cite{liutss3}). It is meaningful to study the minimal {\it $P$-symmetric} periodic solutions of (\ref{1}). So far there are very few papers about it.

In the following, we always suppose $P \in Sp(2n)$ satisfies the $(P)_{k}$ condition.

In this paper, combining the Galerkin approximation procedure (cf. \cite{liutss1,liutss2,liutss3}) with the method with C. Liu and me (cf. \cite{liutss2}), we study the minimal {\it $P$-symmetric} periodic solutions of
(\ref{1}) when the Hamiltonian function $H$ is superquadratic, asymptotically linear and subquadratic respectively.

For $\tau>0$, we define
$$S_{\tau}(H)=\{x\in C^{1}([0, \tau], \R^{2n}): x\neq constant, x \ \text{is a {\it $P$-solution} of }\ (\ref{1})\}.$$

We now state the main results as follows.
\begin{theorem}\label{the:2}
Suppose $P \in Sp(2n)$ satisfies the $(P)_{k}$ condition, and $H$ satisfies the following conditions:
\begin{enumerate}

\item[(H0)] $H\in C^{2}(\R^{2n}, \R)$ with $H(Px) = H(x)$, $\forall x \in \R^{2n}$;

\smallskip
\item[(H1)] $H(x)=\frac{1}{2}(h_{0}x, x)+o(\vert x\vert^{2})$ as $\vert x\vert\to 0$;

\smallskip
\item[(H2)] $H(x)-\frac{1}{2}(h_{0}x, x)\geq 0$, $\forall x\in \R^{2n}$,\\
 where $h_{0}$ is semi-positive definite symmetric matrix with $P^{T}h_{0}P=h_{0}$;

\smallskip
\item[(H3)]  There exist constants $\mu > 2$ and $R_{0} > 0$ such that
\begin{equation*}
0 < \mu H(x) \leq  H^{\prime}(x)\cdot x, \ \ \forall \ \vert x \vert \geq R_{0};
\end{equation*}

\smallskip
\item[(HX1)]  $H^{\prime\prime}(x(t))\geq 0$ for every $x\in S_{\tau}(H)$ and $t \in \R$;

\smallskip
\item[(HX2)]  $\int^{\tau}_{0}H^{\prime\prime}(x(t))dt>0$ for every $x\in S_{\tau}(H)$;

\smallskip
\item[(HX3)] $i_{P}(h_{0})+\nu_{P}(h_{0})\leq \dim\ker_{\R}(P-I)$, where $(i_{P}(h_{0}), \nu_{P}(h_{0}))$ denote the Maslov P-index of $h_{0}$.

\end{enumerate}
Then (\ref{1}) possesses a {\it $P$-solution} $x$ with the minimal {\it $P$-symmetric} period $k\tau$ or $\frac{k\tau}{k+1}$.
\end{theorem}
\begin{remark}
Specially, if $h_{0}=0$, then $i_{P}(h_{0})=0$, $\nu_{P}(h_{0})=\dim\ker_{\R}(P-I)$, $\forall P\in Sp(2n)$. At the moment, (HX3) holds automatically. Our result generalize the corresponding one in \cite{LC2}.
\end{remark}

For the asymptotically linear Hamiltonian systems, we consider the case that the asymptotical matrix may be degenerate and the get the following two theorems:
\begin{theorem}\label{the:8}
Suppose $P \in Sp(2n)$ satisfies the $(P)_{k}$ condition, and $H$ satisfies (H0),(H1), (H2), (HX1), (HX2) and the following conditions:
\begin{enumerate}

\item[(H4)]  There exists constant $a_{1}$, $a_{2}$ and some $s\in (1, \infty)$ such that
\begin{equation*}
\vert H^{\prime\prime}(x)\vert\leq a_{1}\vert x\vert^{s}+a_{2};
\end{equation*}

\item[(H5)]  There exists  semi-positive definite symmetric matrix $h_{\infty}$ with $P^{T}h_{\infty}P=h_{\infty}$ such that
 \begin{equation*}
 H^{\prime}(x)=h_{\infty}x+o(\vert x\vert)\ \  \text{as}\ \  \vert x\vert\to \infty;
 \end{equation*}

\smallskip
\item[(H6)] $h_{\infty}-h_{0}$ is positive definite, $h_{\infty}h_{0}=h_{0}h_{\infty}$, where $h_{0}\in \mathfrak{L}_{s}(\R^{2n})$ is the matrix given in (H1) and (H2);

\smallskip
\item[(HX4)] $i_{P}(h_{\infty})>i_{P}(h_{0})+\nu_{P}(h_{0})$, $i_{P}(h_{0})+\nu_{P}(h_{0})\leq \dim\ker_{\R}(P-I)$, where $(i_{P}(h_{\infty}), \nu_{P}(h_{\infty}))$ denote the Maslov P-index of $h_{\infty}$.

\end{enumerate}
Then (\ref{1}) possesses a {\it $P$-solution} $x$ with the minimal {\it $P$-symmetric} period $k\tau$ or $\frac{k\tau}{k+1}$ provided one of the following cases occurs:
\begin{enumerate}

\item[(1)]$\nu_{P}(h_{\infty})=0$;

\smallskip
\item[(2)]$\nu_{P}(h_{\infty})>0$ and $G_{\infty}(x)=H(x)-\frac{1}{2}(h_{\infty}x, x)$ satisfies
\begin{equation}
\vert G_{\infty}^{\prime}(x)\vert\leq M \ \ \text{for}\ \ x\in\R^{2n}, \ \ G_{\infty}(x)\to +\infty \ \ \text{as}\ \ \vert x\vert\to \infty.
\end{equation}

\end{enumerate}

\end{theorem}

\begin{theorem}\label{the:3}
Suppose $P \in Sp(2n)$ satisfies the $(P)_{k}$ condition, and $H$ satisfies (H0),(H1), (H2),(H4), (H5), (HX1), (HX2) and the following conditions:
\begin{enumerate}

\item[(H7)] $\{x\in \R^{2n}: H^{\prime}(x)=0\}=\{0\}$;

\smallskip
\item[(HX5)]$i_{P}(h_{\infty})+\nu_{P}(h_{\infty})\leq\dim\ker_{\R}(P-I)+1$, $i_{P}(h_{\infty})+\nu_{P}(h_{\infty})\notin [i_{P}(h_{0}), i_{P}(h_{0})+\nu_{P}(h_{0})]$.

\end{enumerate}
Then (\ref{1}) possesses a {\it $P$-solution} $x$ with the minimal {\it $P$-symmetric} period $k\tau$ or $\frac{k\tau}{k+1}$ provided one of the following cases occurs:
\begin{enumerate}

\item[(1)]$\nu_{P}(h_{\infty})=0$;

\smallskip
\item[(2)]$\nu_{P}(h_{\infty})>0$ and $G_{\infty}(x)=H(x)-\frac{1}{2}(h_{\infty}x, x)$ satisfies
\begin{equation}\label{11}
\vert G_{\infty}^{\prime}(x)\vert\leq M \ \ \text{for}\ \ x\in\R^{2n}, \ \ G_{\infty}(x)\to +\infty \ \ \text{as}\ \ \vert x\vert\to \infty.
\end{equation}

\end{enumerate}

\end{theorem}
\begin{remark}
In Theorem \ref{the:3}, we do not need the condition (H6).
\end{remark}

The following theorem studies the minimal {\it $P$-symmetric} periodic solutions of subquadratic Hamiltonian systems with $P$-boundary
\begin{equation}\label{37}
\left\{
\begin{array}{l}
\dot{x}=\lambda JH^{\prime}(x),\ \forall x\in\R^{2n},\ \lambda\in\R,\\
\!x(\tau) = Px(0).
\end{array}\right.
\end{equation}
This is motivated by \cite{BR,FQ2}.

\begin{theorem}\label{the:4}
Suppose $P \in Sp(2n)$ satisfies the $(P)_{k}$ condition, and $H$ satisfies (H0) and
\begin{enumerate}

\item[(H8)]$\vert H^{\prime}(x)\vert\leq M$ for $x\in \R^{2n}$, and $H(x)\to +\infty$ as $\vert x\vert\to\infty$;

\smallskip
\item[(H9)]$H(0)=0$ and $H^{\prime}(x)$, $H(x)>0$ for $x\neq 0$.

\end{enumerate}
Suppose $\tau>0$, (HX1) and (HX2) hold. There exists $\lambda_{\tau}>0$ such that for any $\lambda\geq \lambda_{\tau}$, (\ref{37}) possesses a {\it $P$-solution} $x$  with the minimal {\it $P$-symmetric} period $k\tau$ or $\frac{k\tau}{k+1}$.
\end{theorem}

In order to get the information about the Maslov P-index of the {\it $P$-solution}, we need the relation between the Maslov P-index and Morse index. This has been done in Section 2 by using the Galerkin approximation procedure and the Maslov P-index theory. The main idea comes from \cite{FQ2} and \cite{LC2}.

\bigskip

\section{Maslov P-index and Morse index}
\bigskip

\ \ \ \ Maslov P-index was first studied in \cite{dong} and \cite{LC} independently for any symplectic matrix $P$ with different treatment, it was generalized by C. Liu and the author in \cite{liutss1,liutss2}. And then C. Liu used relative index theory to develop Maslov P-index in \cite{LC2} which is consistent with the definition in \cite{liutss1,liutss2}. In fact, when the symplectic matrix $P=diag \{-I_{n-\kappa}, I_{\kappa}, -I_{n-\kappa}, I_{\kappa}\}$, $0\leq\kappa\in\Z\leq n$, the $(P,\omega)$-index theory and its iteration theory were studied in \cite{donglong} and then be successfully used to study the multiplicity of closed characteristics on partially symmetric convex compact hypersurfaces in $\R^{2n}$. Here we use the notions and results in \cite{LC2,liutss1,liutss2}.

For $\tau>0$, $P\in Sp(2n)$, $\mathfrak{L}_{s}(\R^{2n})$ denotes all symmetric real $2n\times 2n$ matrices. For $B(t)\in C(\R, \mathfrak{L}_{s}(\R^{2n}))$ and satisfies $P^{T}B(t+\tau)P = B(t)$. If $\gamma$ is the fundamental solution of the linear Hamiltonian systems
\begin{equation}\label{2}
\dot y(t)= JB(t)y,\ \ \ y\in \R^{2n}.
\end{equation}
Then the Maslov $P$-index pair of $\gamma$ is defined as a pair of integers
$$(i_{P}, \nu_{P})\equiv (i_{P}(\gamma), \nu_{P}(\gamma))\in \Z\times \{0,1,\cdots,2n\},$$
where $i_{P}$ is the index part and
$$ \nu_{P}=\dim \ker(\gamma(\tau)-P)$$
is the nullity. We also call $(i_{P}, \nu_{P})$ the Maslov P-index of $B(t)$, just as in \cite{LC2,liutss1,liutss2}. If $(\tau, x)$ is a {\it $P$-solution} of (\ref{1}), then the Maslov P-index of the solution $x$ is defined to be the Maslov P-index of $B(t)=H^{\prime\prime}(x(t))$ and denoted by $(i_{P}(x), \nu_{P}(x))$.

Let $S_{k\tau}=\R/(k\tau\Z)$ and $W_{P}=\{z\in W^{1/2,2}(S_{k\tau}, \R^{2n})\mid z(t+\tau)=Pz(t)\}$, it is a closed subspace of $W^{1/2,2}(S_{k\tau}, \R^{2n})$ and is also a Hilbert space with norm $\Vert\cdot\Vert$ and inner product $\langle\cdot, \cdot\rangle$ as in $W^{1/2,2}(S_{k\tau}, \R^{2n})$. Let $\mathfrak{L}_{s}(W_{P})$ and  $\mathfrak{L}_{c}(W_{P})$ denote the space of the bounded selfadjoint linear operator and compact linear operator on $W_{P}$. We define two operators $A$, $B\in \mathfrak{L}_{s}(W_{P})$ by the following bilinear forms:
\begin{equation}\label{3}
\langle Ax, y \rangle = \int_{0}^{\tau}(-J\dot{x}(t), y(t))dt,\ \ \langle Bx, y \rangle = \int_{0}^{\tau}(B(t)x(t), y(t))dt.
\end{equation}

Suppose that $\cdots \leq \lambda_{-j}\leq \cdots \leq \lambda_{-1}<0<\lambda_{1}\leq \cdots \leq \lambda_{j}\leq \cdots$ are all nonzero eigenvalues of the operator $A$ (count with multiplicity), correspondingly, $e_{j}$ is the eigenvector of $\lambda_{j}$ satisfying $\langle e_{j}, e_{i}\rangle=\delta_{ji}$.
We denote the kernel of the operator $A$ by $W_{P}^{0}$ which is exactly the space $\ker_{\R}(P-I)$.
For $m \in \N$, we define the finite dimensional subspace of $W_{P}$ by
\begin{equation*}
W^{m}_{P} = W_{m}^{-}\oplus W_{P}^{0}\oplus W_{m}^{+}
\end{equation*}
with $W_{m}^{-}=\{z\in W_{P}\vert z(t)=\sum_{j=1}^{m}a_{-j}e_{-j}(t), a_{-j}\in\R\}$ and $W_{m}^{+}=\{z\in W_{P}\vert z(t)=\sum_{j=1}^{m}a_{j}e_{j}(t),\\
a_{j}\in\R\}$.

We suppose $P_{m}$ be the orthogonal projections $P_{m}: W_{P} \to W_{P}^{m}$ for $m\in\N\cup \{0\}$. Then $\{P_{m} \mid m=0, 1, 2, \cdots\}$ be the Galerkin approximation sequence respect to $A$.

For $S\in \mathfrak{L}_{s}(W_{P})$, we denote by $M^{*}(S)$ the eigenspaces of $S$
with eigenvalues belonging to $(0, +\infty)$, $\{0\}$ and $(-\infty, 0)$ with  $* = +,0$ and $* = -$, respectively. Similarly, for any $d > 0$, we denote by $M_{d}^{*}(S)$ the $d$-eigenspaces of $S$ with eigenvalues belonging to $[d, +\infty)$, $(-d, d)$ and
$(-\infty, -d]$ with  $* = +,0$ and $* = -$, respectively. We denote $m^{*}(S)=\dim M^{*}(S)$, $m_{d}^{*}(S)=\dim M_{d}^{*}(S)$ and $S^{\sharp} = (S\vert_{Im S})^{-1}$.

The following theorem gives the relationship between the Maslov $P$-index and the Morse index. When $P$ is a symplectic orthogonal matrix, C.Liu in \cite{LC} has got corresponding result. Now we generalize it for any  symplectic matrix $P$. It plays a key role in the proof of the main results.
\begin{theorem}\label{the:1}
Suppose $B(t)\in C(\R, \mathfrak{L}_{s}(\R^{2n}))$ and satisfies $P^{T}B(t+\tau)P = B(t)$ with the Maslov P-index $(i_{P}(B), \nu_{P}(B))$, for any constant $0<d<\frac{1}{4}\Vert (A - B)^{\sharp} \Vert^{-1}$, there exists an $m_{0}>0$ such that for $m\geq m_{0}$, there holds
\begin{equation}\label{35}
\begin{split}
m_{d}^{+}(P_{m}(A - B)P_{m}) &= m+\dim\ker_{\R}(P-I)-i_{P}(B)-\nu_{P}(B),\\
m_{d}^{-}(P_{m}(A - B)P_{m}) &= m+i_{P}(B),\\
m_{d}^{0}(P_{m}(A - B)P_{m}) &= \nu_{P}(B),
\end{split}
\end{equation}
where $B$ be the operator defined by (\ref{3}) corresponding to $B(t)$.
\end{theorem}
\noindent\begin{proof}
Let $x(t) = \gamma_{P}(t)\xi(t)\in W_{P}$, $\xi \in W^{1/2,2}(S_{\tau}, \R^{2n})$, $\gamma_{P}(t)$ is defined in \cite{liutss1,liutss3} is a symplectic path which satisfies $\gamma_{P}(0)=I$ and $\gamma_{P}(\tau)=P$.  Then we have
\begin{equation*}
\begin{split}
\langle Ax, x \rangle &= \int_{0}^{\tau}(-J\dot{x}(t), x(t))dt\\
&= \int_{0}^{\tau}[(-J\dot{\xi}(t), \xi(t)) - (\gamma_{P}(t)^{T}J\dot{\gamma}_{P}(t)\xi(t), \xi(t))]dt\\
&= \int_{0}^{\tau}[(-J\dot{\xi}(t), \xi(t)) - (\bar{B}_{\gamma_{P}}(t)\xi(t), \xi(t))]dt,\\
\langle (A - B)x, x \rangle &= \int_{0}^{\tau}[(-J\dot{x}(t), x(t)) - (B(t)x(t), x(t))]dt\\
&= \int_{0}^{\tau}[(-J\dot{\xi}(t), \xi(t)) - (\gamma_{P}(t)^{T}J\dot{\gamma}_{P}(t)\xi(t), \xi(t))- (\gamma_{P}(t)^{T}B(t)\gamma_{P}(t)\xi(t), \xi(t))]dt\\
&= \int_{0}^{\tau}[(-J\dot{\xi}(t), \xi(t)) - (\widetilde{B}_{\gamma_{P}}(t)\xi(t), \xi(t))]dt,
\end{split}
\end{equation*}
where $\bar{B}_{\gamma_{P}}(t) = \gamma_{P}(t)^{T}J\dot{\gamma}_{P}(t)$, $\widetilde{B}_{\gamma_{P}}(t)=\gamma_{P}(t)^{T}J\dot{\gamma}_{P}(t)+\gamma_{P}(t)^{T}B(t)\gamma_{P}(t)$.
By the definitions of $\gamma_{P}(t)$ and $B(t)$, $\widetilde{B}_{\gamma_{P}}(t)$ and $\bar{B}_{\gamma_{P}}(t)$ are both symmetric matrix functions and $\widetilde{B}_{\gamma_{P}}(0)=\widetilde{B}_{\gamma_{P}}(\tau)$, $\bar{B}_{\gamma_{P}}(0)=\bar{B}_{\gamma_{P}}(\tau)$.
The operators $A$ and $A-B$ defined in $W_{P}$ correspond to the operators $-J\frac{d}{dt}-\bar{B}_{\gamma_{P}}$ and $-J\frac{d}{dt}-\widetilde{B}_{\gamma_{P}}$ defined in $W^{1/2,2}(S_{\tau}, \R^{2n})$. Suppose $\gamma$ is the fundamental solution of $\dot{z}(t)=JB(t)z(t)$.

Consider the following linear Hamiltonian systems
\begin{equation}\label{32}
\dot{z}(t) = J\widetilde{B}_{\gamma_{P}}(t)z(t),\ \ z(t)\in \R^{2n}.
\end{equation}
Suppose $\widetilde{\gamma}(t)$ is the fundamental solution of (\ref{32}). Then by direct computation, we obtain
\begin{equation*}
\widetilde{\gamma}(t) = \gamma_{P}(t)^{-1}\gamma(t)=\gamma_{2}.
\end{equation*}
And similarly, $\gamma_{P}(t)^{-1}$ is the fundamental solution of $\dot{z}(t)=J\bar{B}_{\gamma_{P}}(t)z(t)$.
By Theorem 7.1 in \cite{LO}, there exists an $m^{\ast} > 0$ such that for $m \geq m^{\ast}$ such that
\begin{equation}\label{36}
\begin{split}
m_{d}^{+}(P_{m}(A - B)P_{m}) &= m + i(\bar{B}_{\gamma_{P}})-i(\widetilde{B}_{\gamma_{P}})+\nu(\bar{B}_{\gamma_{P}})- \nu(\widetilde{B}_{\gamma_{P}}),\\
m_{d}^{-}(P_{m}(A - B)P_{m}) &= m -i(\bar{B}_{\gamma_{P}}) + i(\widetilde{B}_{\gamma_{P}}),\\
m_{d}^{0}(P_{m}(A - B)P_{m}) &= \nu(\widetilde{B}_{\gamma_{P}})
\end{split}
\end{equation}
where $\bar{B}_{\gamma_{P}}$ and $\widetilde{B}_{\gamma_{P}}$ be the compact operator defined by (\ref{3}) corresponding to $\bar{B}_{\gamma_{P}}(t)$ and $\widetilde{B}_{\gamma_{P}}(t)$. $(i(\bar{B}_{\gamma_{P}}), \nu(\bar{B}_{\gamma_{P}}))$ and $(i(\widetilde{B}_{\gamma_{P}}), \nu(\widetilde{B}_{\gamma_{P}}))$ is the Maslov-type index of $\bar{B}_{\gamma_{P}}(t)$ and $\widetilde{B}_{\gamma_{P}}(t)$ in \cite{LO}.
Now by Theorem 3.3 in \cite{liutss1}, we have
\begin{equation}\label{33}
i(\bar{B}_{\gamma_{P}})=i_{P}(0)-i(\gamma_{P})-n=-i(\gamma_{P})-n,\ i(\widetilde{B}_{\gamma_{P}})=i_{P}(B)-i(\gamma_{P})-n.
\end{equation}
Note that
\begin{equation}\label{34}
\nu({\bar{B}_{\gamma_{P}}})=\nu(\gamma_{P}(t)^{-1})=\dim\ker_{\R}(P-I),\ \ \nu(\widetilde{B}_{\gamma_{P}}) = \nu(\widetilde{\gamma}) =\nu_{P}(\gamma_{2})=\nu_{P}(\gamma)=\nu_{P}(B).
\end{equation}
Finally we get (\ref{35}) by (\ref{36})-(\ref{34}).
\end{proof}

The following theorem was proved in \cite{LC2} by relative index theory and iteration theory of Maslov P-index.
\begin{theorem}\label{the:7}
Suppose $H\in C^{2}(\R^{2n}, \R)$ and $P \in Sp(2n)$ satisfies the $(P)_{k}$ condition. For $\tau>0$,
let $x_{0}$ be a
{\it $P$-solution} of (\ref{1}). If the Maslov P-index of $x_{0}$ satisfies
$$i_{P}(x_{0}) \leq \dim\ker_{\R}(P-I)+1,$$
and further satisfies (HX1) and (HX2).
Then the minimal {\it $P$-symmetric} period of $x_{0}$ is $k\tau$ or $\frac{k\tau}{k+1}$.
\end{theorem}

In order to estimate the Maslov $P$-index of a critical point of the functional we considered, we need the following result which was proved in \cite{Ghou,Laze,Soli}.

\begin{theorem}\label{the:5}
Let E be a real Hilbert space with orthogonal decomposition $E = X \oplus Y$, where $\dim X < +\infty$. Suppose $f\in C^{2}(E, \R)$ satisfies the (P.S) condition and the following conditions:
\begin{enumerate}

\item[(F1)]  There exist $\rho$ and $\alpha >0$ such that
\begin{equation*}
f(w) \geq \alpha, \ \ \forall w \in \partial B_{\rho}(0) \cap Y.
\end{equation*}

\smallskip
\item[(F2)]  There exist $e\in \partial B_{1}(0) \cap Y$ and $R > \rho$ such that
\begin{equation*}
f(w) < \alpha, \ \ \forall w \in \partial Q.
\end{equation*}
where $Q= (\overline{B_{R}(0)} \cap X) \oplus \{re \mid  0 \leq r \leq R\}$.

\end{enumerate}
Then
\begin{enumerate}
\item  f possesses a critical value $c \geq \alpha$, which is given by
    \begin{equation*}
    c=\inf_{h\in \Lambda}\max_{w\in Q}f(h(w)),
    \end{equation*}
where $\Lambda = \{h\in C(\overline{Q}, E) \mid h=id \ \text{on} \ \partial Q\}$.

 \smallskip
\item If $f^{\prime\prime}(w)$ is Fredholm for $w\in \mathcal{K}_{c}(f) \equiv\{w\in E: f^{\prime}(w)=0, f(w)=c\}$, then there exists an element $w_{0} \in \mathcal{K}_{c}(f)$ such that the negative Morse index $m^{-}(w_{0})$ and nullity $m^{0}(w_{0})$ of $f$ at $w_{0}$ satisfies
    \begin{equation}\label{49}
    m^{-}(w_{0}) \leq \dim X+1 \leq m^{-}(w_{0}) + m^{0}(w_{0}).
    \end{equation}

\end{enumerate}
\end{theorem}

\begin{definition}\cite{Ghou}\label{def:1}
Let $E$ be a $C^{2}$-Riemannian manifold, $D$ is a closed subset of $E$. A family $\mathcal{F}(\alpha)$ is said to be a homological family of dimension $q$ with boundary $D$ if for some nontrival class $\alpha\in H_{q}(E, D)$ the family $\mathcal{F}(\alpha)$ is defined by
$$\mathcal{F}(\alpha)=\{G\subset E: \alpha\  \text{is in the image of} \ i_{\star}: H_{q}(G, D)\to H_{q}(E, D)\},$$
where $i_{\star}$ is the homomorphism induced by the immersion $i: G\to E$.
\end{definition}

\begin{theorem}\cite{Ghou}\label{the:6}
As in the definition \ref{def:1}, for given $E$, $D$ and $\alpha$, let $\mathcal{F}(\alpha)$ be a homological family of dimension $q$ with boundary $D$. Suppose that $f\in C^{2}(E, \R)$ satisfies (P.S) condition. Define
\begin{equation}\label{15}
c\equiv c(f, \mathcal{F}(\alpha))=\inf_{G\in \mathcal{F}(\alpha)}\sup_{w\in G}f(w).
\end{equation}
Suppose that $sup_{w\in D}f(w)<c$ and $f^{\prime}$ is Fredholm on
\begin{equation}\label{16}
\mathcal{K}_{c}=\{x\in E: f^{\prime}(x)=0, f(x)=c\}.
\end{equation}
Then there exists $x\in \mathcal{K}_{c}$ such that the Morse indices $m^{-}(x)$ and $m^{0}(x)$ of the functional $f$ at $x$ satisfy
$$q-m^{0}(x)\leq m^{-}(x)\leq q.$$
\end{theorem}

\bigskip

\section{Superquadratic Hamiltonian systems}
\bigskip

In this section, we study the minimal $P$-symmetric periodic solution of superquadratic Hamiltonian systems with $P$-boundary conditions. In order to prove Theorem \ref{the:2}, we need the following arguments.

For $z\in W_{P}$, we define
\begin{equation}\label{12}
f(z)=\frac{1}{2} \int_{0}^{k\tau}(-J\dot{z}(t), z(t))dt- \int_{0}^{k\tau}H(z)dt=k(\frac{1}{2}\langle Az, z \rangle-\int_{0}^{\tau}H(z)dt).
\end{equation}
It is well known that $f\in C^{2}(W_{P}, \R)$ whenever
\begin{equation}\label{4}
H\in C^{2}(\R^{2n}, \R)\ \ \text{and}\ \ \vert H^{\prime\prime}(x)\vert\leq a_{1}\vert x\vert^{s}+a_{2};
\end{equation}
for some $s\in (1, \infty)$ and all $x\in \R^{2n}$. Looking for solutions of (\ref{1}) is equivalent to looking for critical points of $f$.

\begin{proof}[\bf Proof of Theorem \ref{the:2}] We carry out the proof in several steps.

\bigskip

\noindent{\bf Step 1.} Since the growth condition (\ref{4}) has not been assumed for $H$, we need to truncate the function $H$ at infinite. We follow the method in Rabinowitz's pioneering work \cite{Rab1}.

Let $K>0$ and $\chi\in C^{\infty}(\R, \R)$ such that $\chi(y) \equiv 1$ if $y\leq K$, $\chi(y) \equiv 0$ if $y\geq K+1$, and $\chi^{\prime}(y)<0$ if $y\in (K, K+1)$, where $K$ is free for now. Set
 \begin{equation}\label{72}
H_{K}(z)=\chi(\vert z \vert)H(z)+(1-\chi(\vert z \vert))R_{K}\vert z \vert^{4},
 \end{equation}
where the constant $R_{K}$ satisfies
\begin{equation*}
R_{K}\geq \max_{K\leq \vert z\vert \leq K+1}\frac{H(z)}{\vert z\vert^{4}}.
\end{equation*}
Then $H_{K}\in C^{2}(\R^{2n}, \R)$, and there is $K_{0}>0$ such that for $K\geq K_{0}$, $H_{K}$ satisfies (H1), (H2) and (\ref{4}) with $s=2$. Moreover a straightforward computation shows (H3) hold with $\mu$ replaced by $\nu = \min\{\mu, 4\}$. Integrating this inequality then yields
\begin{equation}\label{28}
H_{K}(z) \geq a_{1}\vert z \vert^{\nu}-a_{2}
\end{equation}
for all $z\in \R^{2n}$, where $a_{1}$, $a_{2}>0$ are independent of $K$.

Let $G_{K}(z)=H_{K}(z)-\frac{1}{2}(h_{0}z, z)$, then by (\ref{28}) it is easy to show that
\begin{equation}\label{5}
G_{K}(z) \geq a_{3}\vert z \vert^{\nu}-a_{4}
\end{equation}
for all $z\in \R^{2n}$, where $a_{3}$, $a_{4}>0$ are independent of $K$.

Finally, we set
\begin{equation}\label{6}
f_{K}(z)=\frac{1}{2} \int_{0}^{k\tau}(-J\dot{z}(t), z(t))dt - \int_{0}^{k\tau}H_{K}(z)dt=\frac{k}{2}\langle Az, z \rangle-\int_{0}^{k\tau}H_{K}(z)dt,\ \ \forall z \in W_{P},
\end{equation}
then $f_{K} \in C^{2}(W_{P}, \R)$.

\medskip
\noindent{\bf Step 2.} For $m>0$, let $f_{K,m}=f_{K}\vert_{W_{P}^{m}}$. We will show that $f_{K,m}$ satisfies the hypotheses of Theorem \ref{the:5}.

By (H1) and (\ref{72}), for any $\epsilon >0$, there is a $M=M(\epsilon, K)>0$ such that
\begin{equation}\label{7}
G_{K}(z) \leq \epsilon\vert z\vert^{2}+M\vert z\vert^{4},\ \forall\ x\in\R^{2n}.
\end{equation}
Let $B_{0}$ be the operator defined by (\ref{3}) corresponding to $h_{0}$, and let
\begin{equation*}
X_{m}=M^{-}(P_{m}(A-B_{0})P_{m})\oplus M^{0}(P_{m}(A-B_{0})P_{m}),\ \ Y_{m}=M^{+}(P_{m}(A-B_{0})P_{m}).
\end{equation*}
For $z\in Y_{m}$, by (\ref{7}) and the fact that $P_{n}B_{0}=B_{0}P_{n}$ for $n\geq 0$, we have
\begin{equation*}
\begin{split}
f_{K, m}(z) &=\frac{k}{2}\langle (A-B_{0})z, z\rangle - \int_{0}^{k\tau}G_{K}(z)dt\\
&\geq \frac{k}{2}\Vert(A-B_{0})^{\sharp}\Vert^{-1}\Vert z\Vert^{2}-(\epsilon \alpha_{2} + M\alpha_{4}\Vert z \Vert^{2})\Vert z \Vert^{2}.\\
\end{split}
\end{equation*}
So there are constant $\rho = \rho(K)>0$ and $\alpha = \alpha(K)>0$, which are sufficiently small and independent of $m$, such that
\begin{equation}
f_{K,m}(z) \geq \alpha, \ \ \forall z\in \partial B_{\rho}(0)\cap Y_{m}.
\end{equation}

Let $e \in \partial B_{1}(0)\cap Y_{m}$ and set
\begin{equation*}
Q_{m}= \{re\mid 0 \leq r \leq r_{1}\}\oplus(B_{r_{1}}(0)\cap X_{m}),
\end{equation*}
where $r_{1}$ is free for the moment. Let $z=z^{-}+z^{0}\in B_{r_{1}}(0)\cap X_{m}$, then
\begin{equation}\label{29}
\begin{split}
f_{K,m}(z+re)&=\frac{k}{2}\langle (A-B_{0})z^{-}, z^{-}\rangle + \frac{k}{2}r^{2}\langle (A-B_{0})e, e\rangle - \int_{0}^{k\tau}G_{K}(z+re)dt\\
&\leq\frac{k}{2}r^{2}\Vert A-B_{0}\Vert - \frac{k}{2}\Vert(A-B_{0})^{\sharp}\Vert^{-1}\Vert z^{-} \Vert^{2}-\int_{0}^{k\tau}G_{K}(z+re)dt.
\end{split}
\end{equation}

If $r=0$, there holds
\begin{equation}\label{31}
f_{K,m}(z+re) \leq - \frac{k}{2}\Vert(A-B_{0})^{\sharp}\Vert^{-1}\Vert z^{-} \Vert^{2}.
\end{equation}

If $r=r_{1}$ or $\Vert z \Vert=r_{1}$, by (\ref{5}), there holds
\begin{equation}\label{8}
\int_{0}^{k\tau}G_{K}(z+re)dt\geq \int_{0}^{\tau}G_{K}(z+re)dt  \geq a_{3}\int_{0}^{k\tau}\vert z+re \vert^{\nu}dt - k\tau a_{4}\geq  a_{5}(\vert z^{0} \vert^{\nu}+r^{\nu})- a_{6}
\end{equation}

Combining (\ref{29}) with (\ref{8}) yields
\begin{equation*}
f_{K,m}(z+re) \leq a_{7}r^{2} - a_{8}\Vert z^{-} \Vert^{2} - a_{5}(\Vert z^{0} \Vert^{\nu}+r^{\nu}) + a_{6}.
\end{equation*}
So we can choose $r_{1}$ large enough which is independent of $K$ and $m$ such that
\begin{equation}
f_{K,m}(z+re) \leq 0, \ \ \forall z \in \partial Q_{m}.
\end{equation}

Now using the same argument as (\cite{liutss2}, Theorem 4.2), we have $f_{K,m}$ has a critical value $c_{K, m}\geq \alpha$, which is given by
\begin{equation}
c_{K, m} = \inf_{g\in \Lambda_{m}}\max_{w\in Q_{m}}f_{K,m}(g(w)),
\end{equation}
where $\Lambda_{m} = \{g\in C(Q_{m}, W_{P}^{m})\mid g=id \ \text{on} \ \partial Q_{m}\}$.
Moreover, there is a critical point $x_{K, m}$ of $f_{K,m}$ which satisfies
\begin{equation}\label{46}
m^{-}(x_{K, m})\leq \dim X_{m} + 1.
\end{equation}

\medskip
\noindent{\bf Step 3.} Since $id\in \Lambda_{m}$, by (\ref{29}) and (H2) we have
\begin{equation}
c_{K, m} \leq \sup_{w\in Q_{m}}f_{K, m}(w) \leq \frac{k}{2}r_{1}^{2}\Vert A-B_{0}\Vert.
\end{equation}
Then in the sense of  subsequence we have
\begin{equation}
c_{K, m} \to c_{K}, \ \ \alpha \leq c_{K} \leq \frac{k}{2}r_{1}^{2}\Vert A-B_{0}\Vert.
\end{equation}

Using the same argument as (4.40)-(4.43) in \cite{liutss2}, we have  that $f_{K}$ satisfies the (P.S)$^{\ast}$ condition on $W_{P}$, i.e., any sequence $\{z_{m}\} \subset W_{P}$ satisfying $z_{m} \in W_{P}^{m}$, $f_{K, m}(z_{m})$ is bounded and $f_{K, m}^{\prime}(z_{m}) \to 0$ possesses a convergent subsequence in $W_{P}$. Hence in the sense of the subsequence we have
\begin{equation}\label{47}
x_{K, m} \to x_{K},\ \ f_{K}(x_{K}) = c_{K},\ \ f_{K}^{\prime}(x_{K})=0.
\end{equation}

By the standard argument as in \cite{liutss2}, $x_{K}$ is a classical nonconstant {\it $P$-solution} of
\begin{equation}\label{9}
\left\{
\begin{array}{l}
\dot{x}=JH_{K}^{\prime}(x),\ \forall x\in\R^{2n},\\
\!x(\tau) = Px(0).
\end{array}\right.
\end{equation}

Indeed, if $x_{K}(t)$ is a constant solution of (\ref{9}), by (H2), then
\begin{equation}\label{56}
f_{K}(x_{K})= \frac{k}{2}\langle Ax_{K}, x_{K}\rangle -\int_{0}^{k\tau}\frac{1}{2}(h_{0}x_{K}, x_{K})dt- \int_{0}^{k\tau}[H_{K}(x_{K})-\frac{1}{2}(h_{0}x_{K}, x_{K})]dt\leq 0.
\end{equation}
This contradicts to $f_{K}(x_{K})=c_{K}\geq \alpha>0$.

And there is a $K_{0}>0$ such that for all $K\geq K_{0}$, $\Vert x_{K} \Vert_{L^{\infty}} <K$. Then $H^{\prime}_{K}(x_{K}) = H^{\prime}(x_{K})$ and $x_{K}$ is a non-constant {\it $P$-solution} of (\ref{9}).
We denote it simply by $x:=x_{K}$.

\medskip
\noindent{\bf Step 4.} Let $B(t)=H_{K}^{\prime\prime}(x(t))$ and $B$ be the operator defined by (\ref{3}) corresponding to $B(t)$. By direct computation, we get
\begin{equation*}
\langle f_{K}^{\prime\prime}(z)w, w \rangle - k\langle (A-B)w, w \rangle = \int_{0}^{k\tau}[H_{K}^{\prime\prime}(x(t))w, w) - (H_{K}^{\prime\prime}(z(t))w, w)]dt, \ \ \forall w \in W_{P}.
\end{equation*}
Then by the continuous of $H_{K}^{\prime\prime}$,
\begin{equation}\label{d9}
\Vert f_{K}^{\prime\prime}(z) - k(A-B) \Vert \to 0 \ \ \ \text{as} \ \ \parallel z - x \parallel \to 0.
\end{equation}
Let $d=\frac{1}{4}\Vert (A-B)^{\sharp} \Vert^{-1}$. By (\ref{d9}), there exists $r_{0}>0$ such that
\begin{equation*}
\Vert  f_{K}^{\prime\prime}(z) - k(A-B) \Vert <\frac{1}{2}d,\ \ \forall z \in V_{r_{0}}=\{z\in W_{P} :\  \Vert z - x\Vert \leq r_{0}\}.
\end{equation*}
Hence for $m$ large enough, there holds
\begin{equation}\label{d6}
\Vert  f^{\prime\prime}_{K, m}(z) -kP_{m}(A-B)P_{m} \Vert <\frac{1}{2}d,\ \ \forall z \in V_{r_{0}} \cap W^{m}_{P}.
\end{equation}
For $z\in V_{r_{0}} \cap W^{m}_{P}$, $\forall w \in M_{d}^{-}(P_{m}(A-B)P_{m})\setminus \{0\}$, from (\ref{d6}) we have
\begin{equation*}
\begin{split}
\langle f^{\prime\prime}_{K, m}(z)w, w \rangle &\leq k\langle P_{m}(A-B)P_{m}w, w \rangle + \Vert  f^{\prime\prime}_{K, m}(z) - kP_{m}(A-B)P_{m} \Vert \cdot \Vert w \Vert^{2}\\
&\leq -d\Vert w \Vert^{2} + \frac{1}{2}d\Vert w \Vert^{2} = -\frac{1}{2}d\Vert w \Vert^{2}<0.
\end{split}
\end{equation*}
Then
\begin{equation}\label{d8}
\dim M^{-}(f^{\prime\prime}_{K, m}(z))\geq \dim M_{d}^{-}(P_{m}(A-B)P_{m}),\ \ \forall z \in V_{r_{0}} \cap W^{m}_{P}.
\end{equation}
Similary to the proof of (\ref{d8}), for large $m$, there holds
\begin{equation}\label{21}
\dim M^{+}(f^{\prime\prime}_{K, m}(z))\geq \dim M_{d}^{-}(P_{m}(A-B)P_{m}),\ \ \forall z \in V_{r_{0}} \cap W^{m}_{P}.
\end{equation}
By (\ref{46}), (\ref{47}), (\ref{d8})  and Theorem \ref{the:1}, for large $m$ we have
\begin{equation*}
\begin{split}
m + i_{P}(h_{0}) + \nu_{P}(h_{0})+1&\geq \dim X_{m} + 1 \geq m^{-}(x_{K, m}) \\
&\geq \dim M_{d}^{-}(P_{m}(A-B)P_{m}) = m + i_{P}(x).
\end{split}
\end{equation*}
Then by (HX3), we have
\begin{equation}\label{10}
i_{P}(x) \leq i_{P}(h_{0}) + \nu_{P}(h_{0})+1\leq \dim\ker_{\R}(P-I)+1.
\end{equation}

Finally, by (\ref{10}), (HX1), (HX2) and Theorem \ref{the:7}, the proof is completed.

\end{proof}

\bigskip

\section{Asymptotically linear Hamiltonian systems}

\bigskip

\begin{proof}[\bf Proof of Theorem \ref{the:3}]Let $W_{P}$, $A$, $P_{m}$ be as in Section 2, and let $f$ be defined by (\ref{12}). Then (H4) implies that $f\in C^{2}(W_{P}, \R)$. Let $B_{0}$ and $B_{\infty}$ be the operator defined by (\ref{3}) corresponding to $h_{0}$ and $h_{\infty}$ respectively.

For $m>0$, let $f_{m}=f\vert_{W_{P}^{m}}$. We carry out the proof in several steps.

\bigskip
\noindent{\bf Step 1.}By (H1), it is easy to prove that
\begin{equation}\label{38}
f(z)=\frac{k}{2}\langle (A-B_{0})z, z\rangle+o(\Vert z\Vert^{2})\ \ \text{as}\ \ z\to 0.
\end{equation}
Let
\begin{equation*}
X_{m}=M^{-}(P_{m}(A-B_{0})P_{m})\oplus M^{0}(P_{m}(A-B_{0})P_{m}),\ \ Y_{m}=M^{+}(P_{m}(A-B_{0})P_{m}).
\end{equation*}
For $z\in Y_{m}$, by (\ref{38}) and the fact that $P_{n}B_{0}=B_{0}P_{n}$ for $n\geq 0$, there exists $\rho>0$ small enough that
\begin{equation}\label{48}
\begin{split}
f_{K, m}(z)&=\frac{k}{2}\langle (A-B_{0})z, z\rangle+o(\Vert z\Vert^{2})\\
&\geq \frac{k}{2}\Vert(A-B_{0})^{\sharp}\Vert^{-1}\Vert z\Vert^{2}+o(\Vert z\Vert^{2})\\
&\geq \alpha=\frac{k}{4}\Vert(A-B_{0})^{\sharp}\Vert^{-1}\Vert \rho\Vert^{2}>0,\ \ \forall \ z\in\partial B_{\rho}(0)\cap Y_{m}.
\end{split}
\end{equation}

\medskip
\noindent{\bf Step 2.}Since $P_{n}B_{\infty}=B_{\infty}P_{n}$ for $n\geq 0$, it is easy to show that there exists $m_{0}>0$ such that
\begin{equation*}
\ker(A-B_{0})\subset W_{P}^{m}.
\end{equation*}
On the other hand, there is $m_{1}>0$ such that for $m\geq m_{1}$
\begin{equation}
\dim\ker(P_{m}(A-B_{\infty})P_{m})\leq \dim\ker(A-B_{\infty}).
\end{equation}
Then there exists $m_{1}\geq m_{0}$ such that for $m\geq m_{1}$,
\begin{equation}
\ker P_{m}(A-B_{\infty})P_{m}=\ker(A-B_{\infty}).
\end{equation}
This implies that
$$Im P_{m}(A-B_{\infty})P_{m}\subset Im(A-B_{\infty}).$$
Then for any $z\in
Im P_{m}(A-B_{\infty})P_{m}$, we have
$$\Vert P_{m}(A-B_{\infty})P_{m}\Vert=\Vert (A-B_{\infty})\Vert\geq \Vert(A-B_{\infty})^{\sharp}\Vert^{-1}\Vert z\Vert.$$
Then for any $0<d\leq \frac{1}{4}\Vert(A-B_{\infty})^{\sharp}\Vert^{-1}$,
\begin{equation}\label{57}
M_{d}^{*}(P_{m}(A-B_{\infty})P_{m})=M^{*}(P_{m}(A-B_{\infty})P_{m}),\ \ \text{where}\ \ * = +,-,0.
\end{equation}
By Theorem \ref{the:1}, there exist $m_{2}\geq m_{1}$ such that for $m\geq m_{2}$,
\begin{equation}\label{39}
\dim M^{-}(P_{m}(A-B_{\infty})P_{m})=m+i_{P}(h_{\infty}).
\end{equation}
Similarly, there exists $m_{3}>0$ such that for $m\geq m_{3}$,
\begin{equation}\label{40}
\begin{split}
\dim M^{-}(P_{m}(A - B_{0})P_{m}) &= m+i_{P}(h_{0}),\\
\dim M^{0}(P_{m}(A - B_{0})P_{m}) &= \nu_{P}(h_{0}),
\end{split}
\end{equation}
Let $m_{4}=\max\{m_{2}, m_{3}\}$. For $m\geq m_{4}$, by (\ref{39}), (\ref{40}) and (HX4) we have
$$\dim M^{-}(P_{m}(A-B_{\infty})P_{m})>\dim X_{m}.$$
It implies that there exists
\begin{equation}\label{42}
y\in  M^{-}(P_{m}(A-B_{\infty})P_{m})\cap Y_{m},\ \ \Vert y\Vert=1.
\end{equation}
By (\ref{42}), we have $(A - B_{\infty})y\in Y_{m}$, $(A - B_{0})y\in Y_{m}$. For any $z=z_{-}+z_{0}\in X_{m}$,
\begin{equation}\label{43}
\langle (B_{\infty}-B_{0})y, z\rangle=-\langle (A-B_{\infty})y, z\rangle+\langle (A-B_{0})y, z\rangle=0.
\end{equation}
By (H6) we have that $B_{\infty}-B_{0}$ is positive definite and
\begin{equation}\label{44}
\langle(B_{\infty}-B_{0})z_{0}, z_{0}\rangle\geq \lambda_{0}\Vert z_{0}\Vert^{2}, \ \ \text{where}\ \ \lambda_{0}>0,
\end{equation}
\begin{equation}\label{41}
[(A-B_{\infty})-(A - B_{0})]^{2}=(B_{0}-B_{\infty})^{2}=(B_{\infty}-B_{0})^{2}=[(A-B_{0})-(A - B_{\infty})]^{2}.
\end{equation}
(\ref{41}) implies that
\begin{equation}
(A-B_{\infty})(A - B_{0})=(A - B_{0})(A - B_{\infty}).
\end{equation}
Hence
\begin{equation}
\begin{split}
0&=\langle (A-B_{\infty})z_{-}, (A-B_{0})z_{0}\rangle=\langle (A-B_{0})(A-B_{\infty})z_{-}, z_{0}\rangle\\
&=\langle (A-B_{\infty})(A - B_{0})z_{-}, z_{0}\rangle=\langle (A - B_{0})z_{-}, (A-B_{\infty})z_{0}\rangle,
\end{split}
\end{equation}
it implies that $\langle z_{-}, (A-B_{\infty})z_{0}\rangle=0$. Hence
\begin{equation}\label{45}
\begin{split}
\langle(B_{\infty}-B_{0})z_{-}, z_{0}\rangle&=\langle z_{-}, (B_{\infty}-B_{0})z_{0}\rangle\\
&=-\langle z_{-}, (A-B_{\infty})z_{0}\rangle+\langle z_{-}, (A-B_{0})z_{0}\rangle=0.
\end{split}
\end{equation}

Set
\begin{equation}
Q_{m}=\{z=ry+z_{-}+z_{0}\in W_{P}^{m}: z_{-}+z_{0}\in X_{m}, \Vert z_{-}+z_{0} \Vert\leq r_{1}, 0\leq r\leq r_{1}\},
\end{equation}
$r_{1}>0$ will be determined later.
For $z=ry+z_{-}+z_{0}\in Q_{m}$, by (\ref{43}), (\ref{44}), (\ref{45}) and (H5), we have
\begin{align*}
f_{m}(z)&=\frac{k}{2}\langle (A-B_{\infty})z, z\rangle-\int_{0}^{k\tau}G_{\infty}(z)dt\\
&=\frac{k}{2}\langle (A-B_{0})z_{-}, z_{-}\rangle+\frac{kr^{2}}{2}\langle (A-B_{\infty})y, y\rangle\\
&\ \ \ \ \ -\frac{k}{2}\langle (B_{\infty}-B_{0})z_{-}, z_{-}\rangle-\frac{k}{2}\langle (B_{\infty}-B_{0})z_{0}, z_{0}\rangle+o(\Vert z\Vert^{2})\\
&\leq -\frac{k}{2}\Vert(A-B_{0})^{\sharp}\Vert^{-1}\Vert z_{-}\Vert^{2}-\frac{kr^{2}}{2}\Vert(A-B_{\infty})^{\sharp}\Vert^{-1}\Vert y\Vert^{2}-\frac{k\lambda_{0}}{2}\Vert z_{0}\Vert^{2}+o(\Vert z\Vert^{2})\\
&\leq -\frac{k}{2}\min\{\Vert(A-B_{0})^{\sharp}\Vert^{-1}, r^{2}\Vert(A-B_{\infty})^{\sharp}\Vert^{-1}, \lambda_{0}\}\Vert z\Vert^{2}+o(\Vert z\Vert^{2}).
\end{align*}
Then taking $r_{1}>$ to be large enough we have
\begin{equation}\label{50}
f_{m}(z)\leq 0,\ \ \forall \ z\in \partial Q_{m}.
\end{equation}

\medskip
\noindent{\bf Step 2.} Using the same arguments as the proof of Lemma 2.1 in \cite{LL} and Lemma 7.1 in \cite{SZ}, we have that $f_{m}$ satisfies (P.S) condition and $f$ satisfies
(P.S)$^{\ast}$ condition either (H5) with $\nu_{P}(h_{\infty})=0$ or the condition (2) in Theorem \ref{the:8}. By (\ref{48}), (\ref{50}) and Theorem \ref{the:5}, $f_{m}$ has a critical value $c_{m}\geq \alpha$, which is given by
\begin{equation}\label{51}
c_{m} = \inf_{g\in \Lambda_{m}}\max_{w\in Q_{m}}f_{m}(g(w)),
\end{equation}
where $\Lambda_{m} = \{g\in C(Q_{m}, W_{P}^{m})\mid g=id \ \text{on} \ \partial Q_{m}\}$.
Moreover, there is a critical point $x_{m}$ of $f_{m}$ which satisfies
\begin{equation}\label{52}
m^{-}(x_{m})\leq \dim X_{m} + 1.
\end{equation}
Since $id\in \Lambda_{m}$, by (\ref{51}) and (H2) we have
\begin{equation}
c_{m} \leq \sup_{w\in Q_{m}}f_{m}(w) \leq \beta =kr_{1}^{2}\Vert A-B_{0}\Vert.
\end{equation}
Then in the sense of  subsequence we have
\begin{equation}\label{53}
c_{m} \to c, \ \ 0<\alpha \leq c \leq \beta.
\end{equation}
Since $f$ satisfies the (P.S)$^{\ast}$ condition on $W_{P}$, hence in the sense of the subsequence we have
\begin{equation}\label{54}
x_{m} \to x,\ \ f(x) = c,\ \ f^{\prime}(x)=0.
\end{equation}

Now using the same arguments as (\ref{56})-(\ref{10}), by (\ref{52})-(\ref{54}) and (HX4), we have that $x$ is a non-constant {\it $P$-solution} of (\ref{1}) with its Maslov P-index $i_{P}(x)$ satisfying
\begin{equation}\label{55}
i_{P}(x) \leq i_{P}(h_{0}) + \nu_{P}(h_{0})+1\leq \dim\ker_{\R}(P-I)+1.
\end{equation}
The proof is completed by (\ref{55}), (HX1), (HX2) and Theorem \ref{the:7}.

\end{proof}

\begin{proof}[\bf Proof of Theorem \ref{the:3}]
\bigskip
{\bf Step 1.} Let $W_{P}$, $A$, $P_{m}$ be as in Section 2, and let $f$ be defined by (\ref{12}). Then (H4) implies that $f\in C^{2}(W_{P}, \R)$. Let $B_{\infty}$ be the operator defined by (\ref{3}) corresponding to $h_{\infty}$.

For $m>0$, let $f_{m}=f\vert_{W_{P}^{m}}$. Using the same arguments as the proof of Lemma 2.1 in \cite{LL} and Lemma 7.1 in \cite{SZ}, we have that $f_{m}$ satisfies (P.S) condition and $f$ satisfies
(P.S)$^{\ast}$ condition either (H5) with $\nu_{P}(h_{\infty})=0$ or the condition (2) in Theorem \ref{the:3}. Let
\begin{equation*}
X_{m}=M^{-}(P_{m}(A-B_{\infty})P_{m})\oplus M^{0}(P_{m}(A-B_{\infty})P_{m}),\ \ Y_{m}=M^{+}(P_{m}(A-B_{\infty})P_{m}).
\end{equation*}
For $z\in Y_{m}$, by (\ref{11}) we have
\begin{equation}\label{17}
\begin{split}
f_{m}(z) &=\frac{k}{2}\langle (A-B_{\infty})z, z\rangle - \int_{0}^{k\tau}G_{\infty}(z)dt\\
&\geq \frac{k}{2}\Vert(A-B_{0})^{\sharp}\Vert^{-1}\Vert z\Vert^{2}- M_{1}\Vert z \Vert^{2}\\
&\geq \alpha=-\frac{k}{2}\Vert(A-B_{0})^{\sharp}\Vert^{-1}M_{1}^{2}.
\end{split}
\end{equation}
For $z=z_{-}+z_{0}\in X_{m}$,  by (\ref{11}) we have
\begin{equation}\label{13}
\begin{split}
f_{m}(z)&=\frac{k}{2}\langle (A-B_{\infty})z_{-}, z_{-}\rangle  - \int_{0}^{k\tau}G_{\infty}(z)dt\\
&\leq  - \frac{k}{2}\Vert(A-B_{0})^{\sharp}\Vert^{-1}\Vert z_{-} \Vert^{2}+M_{1}\Vert z_{-} \Vert-\int_{0}^{k\tau}G_{\infty}(z_{0})dt.
\end{split}
\end{equation}
Since $B_{\infty}P_{n}=P_{n}B_{\infty}$, there exist $m_{1}>0$ such that for $m\geq m_{1}$,
$$\ker P_{m}(A-B_{\infty})P_{m}=\ker (A-B_{\infty}).$$
Then by (\ref{11}),
\begin{equation}\label{14}
\int_{0}^{k\tau}G_{\infty}(z_{0})dt\to +\infty,\ \ \text{as}\ \ \Vert z_{0}\Vert\to\infty.
\end{equation}
By (\ref{13}) and (\ref{14}), there exist $r_{1}>0$ and $\beta <\alpha$ such that
\begin{equation}\label{18}
f_{m}(z)\leq \beta,\ \ \forall \ \ z\in \partial Q_{m},
\end{equation}
where $Q_{m}=\{z\in X_{m}:\Vert z\Vert\leq r_{1}\}$. The constants $\alpha$, $\beta$ and $r_{1}$ in the above are independent of $m$.

\medskip
\noindent{\bf Step 2.} Let $S=Y_{m}$, then $\partial Q_{m}$ and $S$ homologically link
(cf.\cite{chang}). Let $D=\partial Q_{m}$ and $\delta=[Q_{m}]\in H_{l}(W_{P}^{m}, D)$, where $l=\dim X_{m}$. Then $\delta$ is nontrival and  $\mathcal{F}(\delta)$ defined by Definite \ref{def:1} is a homological family of dimension $l$ with boundary $D$. It is well known that $f_{m}^{\prime}$ is Fredholm on $\mathcal{K}_{c_{m}}$ defined by (\ref{15}) and (\ref{16}). By (\ref{17}) and (\ref{18}), we obtain
$$\sup_{w\in D}f_{m}(w)\leq \beta<\alpha\leq c_{m}=c(f_{m}, \mathcal{F}(\delta)).$$
Then by Theorem \ref{the:6}, there exists $x_{m}\in \mathcal{K}_{c_{m}}$ such that the Morse indices $m^{-}(x_{m})$ and $m^{0}(x_{m})$ of $f_{m}$ at $x_{m}$ satisfies
\begin{equation}\label{19}
\dim X_{m}-m^{0}(x_{m})\leq m^{-}(x_{m})\leq \dim X_{m}.
\end{equation}
Since $Q_{m}\in \mathcal{F}(\delta)$, by (\ref{13}) we have
$$c_{m}\leq \sup_{w\in Q_{m}}f_{m}(w)\leq \frac{k}{2}r_{1}^{2}\Vert A-B_{\infty}\Vert+M_{1}r_{1}=M_{2}.$$
Hence in the sense of subsequence we have
\begin{equation*}
c_{m}\to c,\ \ \alpha\leq c\leq M_{2}.
\end{equation*}
Since $f$ satisfies (P.S)$^{\ast}$ condition, in the sense of subsequence,
\begin{equation}
x_{m}\to x_{0},\ \ f(x_{0})=c,\ \ f^{\prime}(x_{0})=0.
\end{equation}
Using the standard arguments we have $x_{0}$ is a classical {\it $P$-solution} of (\ref{1}).
Now using the same arguments as (\ref{d9})-(\ref{d8}), there exists $r_{2}>0$ such that
\begin{equation}\label{20}
\dim M^{\pm}(f^{\prime\prime}_{m}(z))\geq \dim M_{d}^{\pm}(P_{m}(A-B)P_{m}),\ \ \forall z \in \{z\in W_{P}^{m}:\Vert z - x_{0}\Vert \leq r_{2}\},
\end{equation}
where $B$ be the operator defined by (\ref{3}) corresponding to $B(t)=H^{\prime\prime}(x_{0}(t))$.

By (\ref{57}), (\ref{19})-(\ref{20}) and Theorem \ref{the:1}, there exists $m_{2}\geq m_{1}$ such that for $m\geq m_{2}$,
\begin{equation*}
\begin{split}
m + i_{P}(h_{\infty}) + \nu_{P}(h_{\infty})&= \dim X_{m} \geq m^{-}(x_{m}) \\
&\geq \dim M_{d}^{-}(P_{m}(A-B)P_{m}) = m + i_{P}(x_{0})\\
m + i_{P}(h_{\infty}) + \nu_{P}(h_{\infty})&= \dim X_{m} \leq m^{-}(x_{m}) + m^{0}(x_{m})\\
&\leq \dim (M_{d}^{-}(P_{m}(A-B)P_{m})\oplus M_{d}^{0}(P_{m}(A-B)P_{m}))\\
&= m + i_{P}(x_{0})+\nu_{P}(x_{0}).
\end{split}
\end{equation*}
Thus there holds
\begin{equation}\label{22}
i_{P}(h_{\infty})+\nu_{P}(h_{\infty})-\nu_{P}(x_{0})\leq i_{P}(x_{0})\leq i_{P}(h_{\infty})+\nu_{P}(h_{\infty}).
\end{equation}
Combining (\ref{22}) with (HX5)  yields that $x_{0}\neq 0$, or by (H2) we have
\begin{equation}
B(t)=H^{\prime\prime}(x_{0}(t))=h_{0},\ \text{and}\ i_{P}(x_{0})=i_{P}(h_{0}),\ \nu_{P}(x_{0})=\nu_{P}(h_{0}).
\end{equation}
So (\ref{22}) contradicts to (HX5). Further, we have that $x_{0}$ is non-constant by (H7).

Now our conclusion follows from (\ref{22}), (HX5), (HX1), (HX2) and Theorem \ref{the:7}. The proof is complete.

\end{proof}

\bigskip

\section{Subquadratic Hamiltonian systems}
\bigskip

\begin{proof}[\bf Proof of Theorem \ref{the:4}]
\bigskip

Let $W_{P}$, $A$, $P_{m}$ and $W_{P}^{m}$ be defined as in Section 2, let
\begin{equation}\label{23}
g(z)=\lambda \int_{0}^{k\tau}H(z)dt-\frac{k}{2}\langle Az, z \rangle,\ \ \forall z\in W_{P}.
\end{equation}
Set $g_{m}=g\vert_{W_{P}^{m}}$ for $m>0$, it is easy to prove that $g_{m}$ satisfies (P.S) condition and $g$ satisfies (P.S)$^{\ast}$ condition under the condition (H8)(cf.\cite{BR}). Let
\begin{equation*}
X_{m}=P_{m}(M^{+}(A)),\ \ \  Y_{m}=P_{m}(M^{-}(A)\oplus M^{0}(A)).
\end{equation*}
For $z\in X_{m}$, by (H8), (H9) and (\ref{23}),
$$g(z)\leq \lambda M_{1}\Vert z\Vert-\frac{k}{2}\Vert A^{\sharp}\Vert^{-1}\Vert z\Vert^{2}.$$
So there exists $r_{\lambda}>1$ and $Q_{m}=\{z\in W_{P}^{m}: \Vert z\Vert\leq r_{\lambda}\}$ such that
\begin{equation}\label{26}
g(z)\leq 0,\ \ \forall z\in\partial Q_{m}.
\end{equation}
Let $v\in Q_{m}$ with $\Vert z\Vert=1$ and $S_{m}=v+Y_{m}$. For $z=v+z_{-}+z_{0}\in S_{m}$,
\begin{equation}\label{24}
\begin{split}
g(z)&=\lambda\int_{0}^{k\tau}H(z)dt-\frac{k}{2}\langle Az^{-}, z^{-}\rangle-\frac{k}{2}\langle Av, v\rangle  \\
&\leq  \lambda\int_{0}^{k\tau}H(z)dt+\frac{k}{2}\Vert A^{\sharp}\Vert^{-1}\Vert z^{-} \Vert^{2}-\frac{k}{2}\Vert A\Vert.
\end{split}
\end{equation}
Following \cite{BR}, three cases are needed to be considered.
\begin{enumerate}

\item[Case 1]$\Vert z_{-}\Vert^{2}>a_{0}=3\Vert A^{\sharp}\Vert\Vert A\Vert$. Then by (H10) and (\ref{24}),
    $$g(z)\geq \frac{k}{2}(\Vert A^{\sharp}\Vert^{-1}\Vert z^{-} \Vert^{2}-\Vert A\Vert)\geq \Vert A\Vert.$$

\smallskip
\item[Case 2]$\Vert z_{-}\Vert^{2}\leq a_{0}=3\Vert A^{\sharp}\Vert\Vert A\Vert$ and $\Vert z_{0}\Vert>a_{1}$. Then by (H9) and (\ref{24}),
    $$g(z)\geq \lambda k\tau H(z_{0})-\lambda M_{1}\Vert z_{-}+v\Vert+\frac{k}{2}\Vert A^{\sharp}\Vert^{-1}\Vert z^{-} \Vert^{2}-\frac{k}{2}\Vert A\Vert\geq 1$$
    if $\lambda\geq 1$ and $a_{1}$ is large enough.

\smallskip
\item[Case 3]$\Vert z_{-}\Vert^{2}\leq a_{0}$ and $\Vert z_{0}\Vert\leq a_{1}$. Let $S=v+M^{-}(A)\oplus M^{0}(A)$ and $\Omega=\{z\in S:\Vert z_{-} \Vert^{2}\leq a_{0},\Vert z_{0}\Vert\leq a_{1}\}$, then $\Omega$ is convex and weakly compact. Since $\int_{0}^{k\tau}H(z)dt$ is weakly continuous, it achieves its infimum  $\sigma$ on $\Omega$ at $\widehat{z}=v+\widehat{z}_{-}+\widehat{z}_{0}$. By (H10) and the fact $\widehat{z}\ne 0$, we have $\sigma>0$. Therefore,
    $$g(z)\geq\lambda\sigma+\frac{k}{2}\Vert A^{\sharp}\Vert^{-1}\Vert z^{-} \Vert^{2}-\frac{k}{2}\Vert A\Vert\geq 1$$
    if $\lambda\geq \sigma^{-1}(\frac{k}{2}\Vert A\Vert+1)$ and $z\in \Omega$. Hence
    \begin{equation*}
    g(z)\geq 1,\ \ \forall z\in \Omega_{m}=\{z\in S_{m}:\Vert z_{-} \Vert^{2}\leq a_{0},\Vert z_{0}\Vert\leq a_{1}\}.
    \end{equation*}

\end{enumerate}

Combining the three cases, we have the constants
\begin{equation*}
\lambda_{\tau}=\sigma^{-1}(\frac{k}{2}\Vert A\Vert+1)+1,\ \ \alpha=\min\{\Vert A\Vert, 1\}>0
\end{equation*}
such that for $\lambda\geq \lambda_{\tau}$, we have
\begin{equation}\label{25}
g(z)\geq \alpha,\ \ \forall z\in S_{m}.
\end{equation}
Since $\partial Q_{m}$ and $S_{m}$ homologically link, by Theorem II.1.2 and Definition II.1.2 in \cite{chang}, $\partial Q_{m}$ and $S$ homologically link. By (\ref{26}) and (\ref{25}), using the same argument as step 2 in the proof of Theorem \ref{the:3}, there is a classical {\it $P$-solution} $x_{0}$ of (\ref{37}) such that
\begin{equation}\label{30}
i_{P}(x_{0})\leq \dim\ker_{\R}(P-I)
\end{equation}
\begin{equation}\label{27}
g(x_{0})=c\geq \alpha>0,\ \ g^{\prime}(x_{0})=0.
\end{equation}
By (H9) and (\ref{27}), $x_{0}$ is non-constant. By (\ref{30}), (HX1), (HX2) and Theorem \ref{the:7}, we complete the proof.
\end{proof}


\begin{thebibliography}{10}



\bibitem{HE}
H.~Amann and E.~Zehnder, \emph{Nontrivial solutions for a class of nonresonance problems and applications to nonlinear differential equations}, Ann. Sc. Super. Pisa, Cl. Sci. Serie IV, VII(4) (1980), 539--603.

\bibitem{BR}
V.~Benci and P.~Rabinowitz, \emph{Critical point theorems for indefinite functionals}, Inv. Math. \textbf{52} (1979), 241--273.

\bibitem{CAMR}
A.~Chenciner and R.~Montgomery, \emph{A remarkable periodic solution of the three body problem in the
case of equal masses}, Ann. of Math. \textbf{152} (2000), no.3, 881--901.

\bibitem{chang}
K. C.~Chang, \emph{Infinite dimensional Morse theory and multiple solution problems}, in Progress in Nonlinear Differential Equations and Their Applications, Vol 6 (1993).


\bibitem{CLL}
K. C.~Chang, J.~Liu and M. Liu, \emph{Nontrivial periodic solutions for strong resonance Hamiltonian systems}, Ann. Inst. H. Poincare Anal. Non. lineaire \textbf{14}(1) (1997), 103--117.

\bibitem{dong1}
Y.~Dong, \emph{Maslov type index theory for linear Hamiltonian systems with Bolza boundary value conditions and multiple solutions for nonlinear Hamiltonian systems},  Pacific J. Math. {\bf 221}: 2 (2005), 253--280.


\bibitem{dong}
Y.~Dong, \emph{$P$-index theory for linear Hamiltonian systems and multiple solutions for nonlinear Hamiltonian systems},   Nonlinearity {\bf 19}: 6 (2006), 1275--1294.


\bibitem{donglong}
Y.~Dong and Y.~Long, Closed characteristics on partially symmetric compact convex hypersurfaces in $\R^{2n}$, J. Diff. Equa. {\bf 196} (2004), 226-248.

\bibitem{FDTS}
D.~Ferrario and S.~Terracini, \emph{On the existence of collisionless equivariant minimizers for the classical $n$-body problem}, Invent. Math. \textbf{155}: 2 (2004), 305--362.


\bibitem{FQ1}
G.~Fei and Q.~Qiu, \emph{Periodic solutions of asymptotically linear Hamiltonian systems}, Chin. Ann. Math.,  \textbf{18B}(3) (1997), 359--372.

\bibitem{FQ2}
G.~Fei and Q.~Qiu, \emph{Minimal period solutions of nonlinear Hamiltonian systems}, Nonlinear Analysis, Theory, Method  Applications  \textbf{27}(7) (1996), 821--839.

\bibitem{Ghou}
N.~Ghoussoub, \emph{Location, multiplicity and Morse indices of min-max critical points}, J. Reine Angew Math. \textbf{417} (1991), 27--76.

\bibitem{HuSun2}
X.~Hu and S.~Sun, \emph{Morse index and the stability of closed geodesics}, Sci. China Math. \textbf{53}(5) (2010), 1207--1212.

\bibitem{HuSun}
X.~Hu and S.~Sun, \emph{Index and stability of symmetric periodic orbits in Hamiltonian systems with application to figure-eight orbit}, Comm. Math. Phys. \textbf{290} (2009), 737--777.

\bibitem{HuSun1}
X.~Hu and S.~Sun, \emph{Stability of relative equilibria and Morse index of central configurations}, C. R. Acad. Sci. Paris \textbf{347} (2009), 1309--1312.

\bibitem{HuWang}
X.~Hu and P.~Wang, \emph{Conditional Fredholm determinant for the $S$-periodic orbits in Hamiltonian systems}, Joural of Functional Analysis \textbf{261} (2011), 3247--3278.

\bibitem{Laze}
A.~Lazer and S.~Solomini, \emph{Nontrivial solution of operator equations and Morse indices of critical points of min-max type}, Nonlinear Anal. \textbf{12} (1988), 761--775.


\bibitem{LL}
S.~Li and J.~Liu, \emph{Morse theory and asymptotically linear Hamiltonian systems}, J. Diff. Equa. \textbf{78} (1989), 53--73.

\bibitem{LC}
C.~Liu, \emph{Maslov P-index theory for a symplectic path with applications}, Chin. Ann. Math. \textbf{4} (2006), 441--458.

\bibitem{LC5}
C.~Liu, \emph{Periodic solutions of asymptotically linear delay differential systems via Hamiltonian systems}, J. Differential Equations \textbf{252} (2012), 5712--5734.


\bibitem{LC2}
C.~Liu, \emph{Relative index theories and applications}, 2015.

\bibitem{liutss1}
C.~Liu and S.~Tang, \emph{Maslov $(P, \omega)$-index theory for symplectic paths}, Advanced Nonlinear Studies \textbf{15} (2015), 963--990.

\bibitem{liutss2}
C.~Liu and S.~Tang, \emph{Iteration inequalities of the Maslov $P$-index theory with applications}, Nonlinear Analysis \textbf{127} (2015), 215--234.

\bibitem{liutss3}
C.~Liu and S.~Tang, \emph{Subharmonic P-solutions of first order Hamiltonian systems}, submitted.

\bibitem{LO}
Y.~Long, \emph{Index theory for symplectic paths with application}, Progress in Mathematics, Vol. 207, Birkh$\ddot{a}$user Verlag, 2002.

\bibitem{Rab1}
P. H.~Rabinowitz, \emph{Periodic solutions of Hamiltonian systmes}, Comm. Pure Appl. Math. \textbf{31} (1978), 157--184.


\bibitem{Soli}
S.~Solimini, \emph{Morse index estimates in min-max theorems}, Manuscription Math. \textbf{63} (1989), 421--453.


\bibitem{SZ}
A.~Szulkin, \emph{Cohomology and Morse theory for strongly indefinite functionals}, Math. Z. \textbf{209} (1992), 375--418.









\end{thebibliography}
\end{document}